\documentclass[12pt]{article}

\def\bfX{{\mathbf X}}
\def\wXc{[{\bfX(c)}]}
\def\wX{[{\bfX}]}
\def\bfs{{\mathbf s}}
\usepackage{amsmath,mathrsfs}
\def\MR{MR}

\usepackage{shapepar,bm}
\usepackage[dvips,arrow,matrix,ps,color,line]{xy}
\usepackage{theorem,amsfonts,amssymb,amscd}
\usepackage{geometry}
\usepackage{makeidx}

\providecommand{\bysame}{\leavevmode\hbox to3em{\hrulefill}\thinspace}
\providecommand{\MR}{\relax\ifhmode\unskip\space\fi MR }

\providecommand{\href}[2]{#2}

\usepackage{graphicx}
\usepackage{graphics}

\usepackage{epstopdf}

\usepackage{amssymb,graphics,hyperref}

\hypersetup{
    colorlinks = true,
    linkcolor = blue,
    anchorcolor = red,
   citecolor = blue,
    filecolor = red,
    pagecolor = red,
    urlcolor = blue}

\newcommand{\mathsym}[1]{{}}

\usepackage[usenames]{color}
\definecolor{MyLightMagenta}{cmyk}{0.1,0.8,0,0.1}
\definecolor{MyDarkBlue}{rgb}{0.1,0,0.3}

\makeindex

\def\ovsig{\overline{\sigma}}

\def\bfc{{\mathbf c}}

\def\bfw{{\mathbf w}}

\def\d{\partial}

\def\mod{{\mathrm{mod}}}

\def\ZZ{\mathbb Z}

\def\CC{\mathbb C}

\def\Ecal{{\mathcal E}}

\def\QQ{\mathbb Q}

\def\cocoa{{\hbox{\rm C\kern-.13em o\kern-.07em C\kern-.13em o\kern-.15em A}}}

\def\Dcal{\mathcal D}

\def\Res{{\rm Res}}

\def\ttp{{\tt p}}

\def\End{\mathrm{End}}

\def\blamb{{\bm \lambda}}

\def\Pcal{{\mathcal P}}

\def\w2M{\bigwedge^2M}

\def\w{\wedge }
\def\bw{\bigwedge }

\protect
\protect
\protect
\protect

\protect
\protect

\def\sra{\rightarrow}
\def\lra{\longrightarrow}

\def\proof{\noindent{\bf Proof.}\,\,}
\def\qed{{\hfill\vrule height4pt width4pt depth0pt}\medskip}
\def\be{\begin{equation}}
\def\ee{\end{equation}}
\def\bclm{\begin{claim}}
\def\eclm{\end{claim}}
\def\beqn{\begin{eqnarray}}
\def\eeqn{\end{eqnarray}}
\def\beqn*{\begin{eqnarray*}}
\def\eeqn*{\end{eqnarray*}}

\theoremstyle{change}
\theorembodyfont{\rmfamily}

\newtheorem{claim}{}[section]

\global
\global

\def\no@breaks#1{{\def\\{ \ignorespaces}#1}}    

  \makeatletter
\def\cleardoublepage{\clearpage\if@twoside \ifodd\c@page\else
\hbox{} \thispagestyle{empty}
\newpage
\if@twocolumn\hbox{}\newpage\fi\fi\fi} \makeatother

\usepackage{eso-pic,graphicx}
\makeatletter
\newcommand\BackgroundPicture[2]{%
  \setlength{\unitlength}{1pt}%
  default \put(0,\strip@pt\paperheight){%
  \parbox[t][\paperheight]{\paperwidth}{%
    \vfill
     \centering \includegraphics[angle=#2, width=15cm, height=15cm,  bb=0 0 150 150]{#1}
    \vfill
}}} %
\makeatother

\date{}

\title{Universal Decomposition Algebras Represent Endomorphisms}
\author{{\sc ommolbanin behzad and abbas nasrollah nejad}}

\begin{document}
\footnotetext{
	{\em Keywords and Phrases:} Hasse-Schmidt Derivation, Exterior Algebra, Universal Decomposition Algebra, Giambelli's formula, Bosonic and Fermionic Representations by Date-Jimbo-Kashiwara-Miwa.\\
}
\maketitle

\begin{abstract} \noindent The goal of this paper is to supply an explicit description  of the universal decomposition algebra of the generic polynomial of degree $n$ into the product of two monic polynomials, one of degree $r$,  as a representation of  Lie algebras of $n\times n$ matrices with polynomial entries. This is related with the bosonic vertex representation of the Lie algebra $gl_\infty$ due to Date, Jimbo, Kashiwara and Miwa.

\end{abstract}

\section*{Introduction} 
\claim{\bf The Goal.} Let $X^n(c)$ denote the generic monic polynomial of degree $n\geq 1$
$$
X^n-c_1X^{n-1}+\cdots+(-1)^nc_n\in B_0(c)[X],
$$
where $B_0(c):=\QQ[c_1,\ldots,c_n]$ is the polynomial ring in the $n$ 
indeterminates $(c_1,\ldots,$ $ c_n)$.
The main purpose of this paper is {\em to provide an explicit description of the 
universal decomposition algebra $B_{r,n}(c)$ of $X^n(c)$, $ 1\leq r\leq  n$, 
into the product of two monic polynomials, one of degree $r$, as  a 
representation of the Lie algebra $gl_n(B_0(c))$ of the $B_0(c)$-valued $n
\times n$ matrices}. Our description is achieved by computing the action on 
$B_{r,n}(c)$ of the generating function of the  basis elements of  $gl_n(B_0(c))$ 
formed by the generic monic polynomials $X^i(c)\in V(c):=B_0(c)[X]$ of 
degree  $0\leq i\leq n-1$, combining techniques and results by Laksov and Thorup on one hand and the notion of Schubert derivation on exterior algebras introduced by Gatto in \cite{SCHSD}. Putting all the $c_i=0$ one obtains, as a 
particular case,  \cite[Theorem 4.7]{gln} by Gatto and Salehyan,
describing the singular cohomology $B_{r,n}$ of the complex Grassmannian 
$G(r,n)$ as a module over the Lie algebra of $\QQ$-valued $n\times n$ 
matrices. The latter, as pointed out by those authors, can be seen as a finite dimensional counterpart of the celebrated 
vertex representation of the Lie algebra $gl_n(\CC)$ due to Date, Jimbo, 
Kashiwara and Miwa (DJKM) \cite{DJKM01,jimbomiwa}, which in turn  can be recovered as the limit for $r,n\sra \infty$, as shown e.g. in \cite{SDIWP}. 

\claim{\bf Relationship with related literature.} It is well  known that the singular cohomology of the complex Grassmannian $G(r,n)$ is the universal decomposition algebra $B_{r,n}$ of the polynomial $X^n$ as the product of two monic polynomials, one of degree $r$ (see \cite{LakTh02}).  Thus, our main Theorem \ref{thm44} may be understood as  the extension of \cite[Theorem 4.7]{gln}, describing the cohomology of Grassmannians  as a representation of  the Lie algebra of the $n\times n$ square $\QQ$-valued matrices, to the cohomology of Grassmann bundles  $G(r,E)\sra X$ over a smooth variety of dimension at least $n$, once the coefficients of the generic polynomial $X^n(c)$  are interpreted as the Chern classes of $E\sra X$. We have chosen, however, to adopt in our exposition an  entirely algebraic language, to keep at a minimum the  prerequisites required to the  readership. 

The  approach we propose, consists in blending and, in a sense, reconciling two alternatives ways to look at Schubert Calculus via exterior powers due, on one hand,  to Laksov and Thorup \cite{LakTh01,LakTh02,LakTh06} (see also \cite{LakEq1,LakEq2} for the equivariant version) and, on the other hand, to Gatto and collaborators who rely on the notion of Hasse-Schmidt derivations on a exterior algebra \cite{SCHSD,SCGA,ESC, HSDGA}. The special case relevant for our applications  is the Schubert derivation,  whose kinship with the Boson-Fermion correspondence and the vertex operators arising in the representation theory  of Heisenberg algebras has been highlighted in \cite{HSDGA,pluckercone,gln,SDIWP}. The rest of this introduction shall be devoted to describe more precisely our main results and methods while the last part will be concerned with the content of the paper, section by section.
\claim{} The universal decomposition algebra $B_{r,n}(c)$ of $X^n(c)$ is a quotient of $B_r(c):=B_0(c)[e_1,\ldots,e_r]$. Let $\pi_{r,n}:B_r(c)\sra B_{r,n}(c)$ be the canonical projection and let 
$$E_r(z)=1-e_1z+\cdots+(-1)^re_rz^r.
$$ Let $V(c):=B_0(c)[X]$ and $V_n(c)=V(c)/(X^n(c))$. The first observation is that $B_{r,n}(c)$ is a $gl(V_n(c))$-module. The reason is that by \cite[The main Theorem 0.1]{LakTh01}, there is a $B_0(c)$-module isomorphisms 
\be
B_{r,n}(c)\sra \bw^rV_n(c),
\ee
mapping certain deformations of Schur polynomials $\Delta_\blamb(H_{r,n}
(c))$, forming a basis of $B_{r,n}(c)$ parametrised by partitions whose 
Young diagram is contained in a $r\times (n-r)$ rectangle, to the basis 
element $\wXc^r_\blamb:=X^{r-1+\lambda_1}(c)\w\cdots \w X^{\lambda_r}(c)$ 
(Cf. \ref{thm:lakth}).  
As the exterior $\bw V_n(c)$ represents canonically the $B_0(c)$ 
endomorphisms of $V_n(c)$ (Cf. Section~\ref{sec2:glnwr}) and  the exterior 
powers are irreducible sub-representations, it follows that $\bw^rV_n(c)$ is 
itself  a $gl(V_n(c))$-module for all $r\geq 0$. Thus $B_{r,n}(c)$ gets 
a $gl(V_n(c))$-module structure by transporting that of $\bw^rV_n(c)$ via 
the Laksov-Thorup isomorphism. The core of our paper consists in providing 
the explicit description of  the $B_{r,n}(c)$ representation of $gl(V_n(c))$ 
by exploiting the notion of Schubert derivation introduced in \cite{SCHSD}, 
extensively treated in \cite{HSDGA} and used in \cite{SDIWP} to supply an 
alternative deduction of the generating function accounting for the bosonic 
representation of $gl_\infty$ after DJKM. 
Let $(\d^i(s))_{0\leq i\leq n-1}$ 
be the basis of $V_n(c)^*$, dual of $(X^i(c))_{0\leq i\leq n-1}$. Thus
$$
gl(V_n(c))=\bigoplus_{i,j=0}^{n-1}B_0(c)\cdot X^i(c)\otimes \d^j(s),
$$
and we consider the generating function of the basis of $gl(V_n(c))$, 
$$
\Ecal_n(z,w^{-1})=\displaystyle{\sum_{i,j=0}^{n-1}}X^i(c)\otimes \d^j(s)z^iw^{-j}.
$$ 
Our main result is the description of the action of $\Ecal_n(z,w^{-1})$ on any basis element of $B_{r,n}(c)$.

\medskip
\noindent
{\bf Theorem \ref{thm44}}
\begin{eqnarray*}
	&&\Ecal_n(z,w^{-1})(\Delta_\blamb(H_{r,n}(c)))\cr\cr
	&=&\pi_{r,n}{c(z)z^{r-1} \over c(w)E_r(z)}\begin{vmatrix}w^{-r+1-\lambda_1}(c)&w^{-r+2-\lambda_2}(c)&\cdots &w^{-\lambda_r}(c)\cr\cr
	\ovsig_-(z)h_{\lambda_1+1}(c)&\ovsig_-(z)h_{\lambda_2}(c)&	\ldots&\ovsig_-(z)h_{\lambda_r-r+2}(c)\cr\cr
	\vdots&\vdots&\ddots&\vdots\cr
	\ovsig_-(z)h_{\lambda_1+r-1}(c)&\ovsig_-(z)h_{\lambda_1+r-2}(c)&\ldots&\ovsig_-(z)h_{\lambda_r}(c)\end{vmatrix}
\end{eqnarray*}
where $\sum_{j\in\ZZ}h_j(c)z^j=\displaystyle{1\over E_r(z)}\in B_r(c)[[z]]$ and $\ovsig_-(z)$ (which is the differential part of a vertex operator) acts on $h_j(c)$ according to the rule $\ovsig_-(z)h_j(c)=h_j(c)-h_{j-1}(c)z^{-1}$ and $c(z)=1-c_1z+\cdots+(-1)^nc_nz^n$.
For instance if $n=4$ and $r=2$, we show that
\begin{eqnarray*}
 \Ecal_{4}(z,w^{-1})\Delta_{(2,1)}(H_{2,4}(c))&=&{h_2(c)\over w} + h_1(c)h_2(c){z\over w}+ h_2(c)^2 {z^2\over w} -{1\over w^3} \cr\cr &+& (h_1(c)^2-h_2(c)){z^2\over w^3}
+ h_1(c) h_2(c){z^3\over w^3}\label{eq4:eq25}
\cr\cr\cr
   &=&h_2(c)\,{1\over w}+\Delta_{(2,1)}(H_{2,4}(c)){z\over w}+\Delta_{(2,2)}(H_{2,4}(c)){z^2\over w}\cr\cr
   &-&{1\over w^3}+\Delta_{(1,1)}(H_{2,4}(c)){z^2\over w^3}+ \Delta_{(2,1)}(H_{2,4}(c)){z^3\over w^3}.
\end{eqnarray*}

Theorem \ref{thm44} is obtained from \ref{thm.thm210}, computing the $gl(V(c))$-structure of $B_r(c)$ by checking that the result can be harmless projected on $B_{r,n}(c)$. 

\claim{\bf Organisation of the paper.} The paper as a whole is organised as 
follows. Section \ref{sec1:sec1} is devoted to preliminaries and notation. 
This section is quite necessary because in spite there are very precise 
reference to Hasse-Schmidt derivations of exterior algebras (see also 
\cite{GSCH} in addition to the previously cited ones), most of the terminology is not 
standard yet,  although has already been mentioned by other authors like 
e.g. in \cite{Nigran,tango} or, more recently,  in \cite[p.~116]{CoDuGu}. In 
this section we also recall the theory of Laksov and Thorup together with 
their definition of residue of an $r$-tuple of Laurent polynomials, as in 
\cite{LakTh01}. Section \ref{sec2:sec2} begins by recalling the definition 
of the Schubert derivation through the canonical $\bw V(c)$ representation 
of the Lie algebra $gl(V(c))$ and  in Theorem  \ref{thmgmblt} we recall the 
Laksov and Thorup version of the Giambelli formula holding in $B_r(c)$. In 
Section \ref{sec3:sec3} we describe the $gl(V(c)$-module structure of 
$B_r(c)=B_0(c)[e_1,\ldots,e_r]$. The outputs are Theorem \ref{thm.thm210} 
and  Corollary \ref{cor3:cor22}, whose proofs, far from being a trivial 
repetition of the main theorems of \cite{gln}, need to be suitably adapted 
and require the overpassing of some subtleties. One of the main tools is an 
easy but powerful formalism to describe the contraction of an exterior 
monomial against a linear form, for which we do not know any explicit 
reference to the literature, described in details and with examples in 
Section~\ref{maintool}. The paper concludes with Section \ref{sec4:sec4} 
devoted to our main result of describing through suitable generating 
functions the universal decomposition algebra $B_{r,n}(c)$ as a 
representation of the Lie algebra $gl(V_n(c))$. Here the main issue is to 
check that the generating function $\Ecal_n(z,w^{-1})$ for $n=\infty$ 
preserves the ideal of relations passing from $B_r(c)$ to $B_{r,n}(c)$. The 
other consists in getting rid in a efficient way of the combinatorics needed 
to prove Theorem~\ref{thm44}, which we then illustrates with an explicit 
Example for special low values of $r$ and $n$. This enables the readers to 
check that the  formula  works by computing the image of a basis element of 
$B_{r,n}(c)$ through a specific basis element of $gl(V_n(c))$. As we have said in case $c_i=0$ one gets back the main theorems of \cite{gln} and the case $r,n=\infty$ does correspond, as explained in \cite{gln,SDIWP}, with the DJKM representation of the Lie algebra $gl_\infty(\CC)$ on a polynomial ring $B$ (the bosonic Fock space) of polynomials in infinitely many indeterminates (as pedagogically explained  e.g. in \cite[Section 5.4]{KacRaRoz}).

\section{Preliminaries and Notations}\label{sec1:sec1}
\claim{} Let $B_0(c):=\QQ[c_1,\ldots,c_n]$ be the polynomial ring in the $n$ 
indeterminates $(c_1,\ldots,$ $ c_n)$ and  $V(c):=B_0(c)[X]$ the polynomial 
ring with coefficients in $B_0(c)$ in the indeterminate $X$. Note that  $V(c) $ is a free  $B_0(c)$-module with basis either $\bfX:=(X^i)_{i\geq 0}$ or  $
\bfX(c):=(X^j(c))_{j\geq 0}$, where 
$
X^j(c)=X^j+\sum_{i=1}^j(-1)^ic_iX^{j-i}.
$
By convention, we put $X^j=0$ if $j<0$ and $c_j=0$ if $j>n$.
Setting 
\be
\bfc(z):=1-c_1z+c_2z^2+\cdots+(-1)^rc_rz^r.
\ee
Let $V(c)^*:=\bigoplus B_0(c)\d^j$ be the restricted dual of $V(c)$, where
$$
\d^j:={1\over j!}\left.{d^j\over dX^j}\right|_{X=0},
$$
such that $(\d^j(X^i)=\delta_j^i)$. It follows that
\be
gl(V(c)):=V(c)\otimes V(c)^*=\bigoplus_{i,j\geq 0}B_0(c)\cdot X^i\otimes \d^j.\label{eq1:xidi}
\ee
\claim{} We denote by $\bw V(c)$ the exterior algebra of $V(c)$. Notice that $\bw^0V(c)=B_0(c)$. For $r\geq 1$ we shall use the following two bases $(\wX^r_\blamb)_{\blamb\in\Pcal_r}$ and $(\wXc^r_\blamb)_{\blamb\in\Pcal_r}$, where $\Pcal_r$ denotes the set of all partitions of length at most $r$, and
$$
\wX^r_\blamb:=X^{r-1+\lambda_1}\w X^{r-2+\lambda_2}\w\cdots\w X^{\lambda_r},
$$
and
$$
\wXc^r_\blamb:=X^{r-1+\lambda_1}(c)\w X^{r-2+\lambda_2}(c)\w\cdots\w X^{\lambda_r}(c).
$$
\claim{\bf Hasse-Schmidt Derivations on $\bw V(c)$} Recall from \cite{HSDGA}
that a Hasse-Schmidt derivation n $\bw V(c)$ is a map $\Dcal(z): \bw V(c)
\sra \bw V(c)[[z]]$ such that $\Dcal(z)(u\w v)=\Dcal(z)u\w \Dcal(z)v$ for 
all $u,v\in \bw V(c)$. By \cite{HSDGA, GSCH} any endomorphism $f$ of $V(c)$ 
uniquely determines a HS-derivation, whose restriction to $V(c)=\bw^1V(c)$ is precisely $f$.
\claim{\bf The universal decomposition algebra.} \label{sec1:uda}
Let $(e_1,\ldots,e_r)$ be indeterminates over $B_0(c)$. Let $B_r(c):=B_0(c)[e_1,\ldots,e_r]$. In the ring $B_r(c)[z]$ consider the generic polynomial  $E_r(z):=1-e_1z+\cdots+(-1)^re_rz^r$ and $H_r(z)$ implicitly defined by 
\be
H_r(z)=\sum_{j\in\ZZ}h_jz^j={1\over E_r(z)}\in B_r(c)[[z]].\label{eq:Hrz}
\ee
In particular $h_j=0$ if $j<0$ and  $h_0=1$.
Let $H_r(c):=(h_0(c),h_1(c),\ldots)$ be the sequence defined by the equality
\be
\sum_{j\in\ZZ}h_j(c)z^j=\bfc(z)H_r(z)={\bfc(z)\over E_r(z)}.
\ee
In particular $h_j(c)=0$ if $j<0$, $h_0(c)=1$ and $h_1(c)=h_1-c_1$, 
$h_2(c)=h_2-c_1h_1+c_2$, \ldots, $h_k(c)=h_k-\sum_{j\geq 0}c_jh_{k-j}$, with 
the convention that $c_j=0$ if $j\geq n-1$. We also need the notation
\be
\ttp_r(X)=X^rE_r\left(X^{-1}\right)\label{eq1:ttprX}.
\ee
Let 
\be
B_{r,n}(c):={B_r(c)\over (h_{n-r+1}(c),\ldots, h_n(c))},\label{eq2:brnc}
\ee
 and $\pi_{r,n}:B_r(c)\sra B_{r,n}(c)$, where $B_r(c):=B_0(c)[e_1,\ldots,e_r]$.
\bclm{\bf Proposition.} {\em The $B_0(c)$-algebra $B_{r,n}(c)$ is the universal decomposition algebra of the generic monic polynomial $X^n(c)$.}
\eclm
\proof The claim means that $X^n(c)$ decomposes in $B_{r,n}(c)[X]$ as the product of two monic polynomials, one of degree $r$. To see this, firts notice that in $B_{r,n}(c)$ we have the following decomposition
$$
\bfc(z)=E_r(z)(1+h_1(c)z+h_2(c)z^2+\cdots+h_{n-r}(c)z^{n-r}),
$$
using the variable $z=X^{-1}$ we have
$$
{X^n(c)\over \ttp_r(X)}=X^{n-r}(1+h_1(c)X^{-1}+\cdots+h_{n-r}(c)X^{n-r})=X^{n-r}+h_1(c)X^{n-r-1}+\cdots+h_{n-r}(c).
$$
We claim that
\be
\ttp_r(X)(X^{n-r}+h_1(c)X^{n-r-1}+\cdots+h_{n-r}(c)),\label{eq1:undec}
\ee
is the universal decomposition of $X^n(c)$ in $B_{r,n}(c)$.
 Indeed,  if $A$ is any $B_0(c)$--algebra where $X^n(c)=p(X)q(X)$ for some  monic polynomials 
$$
p(X)=X^r-e_1(p)X^{r-1}+\cdots+(-1)^re_r(p),\qquad e_i(p)\in A,
$$
and 
$$
q(X)=X^{n-r}+h_1(q)X^{n-r-1}+\cdots+h_{n-r}(q),\qquad h_i(q)\in A,
$$
of degree $r$ and $n-r$ respectively, then the unique $B_0(c)$-algebra 
homomorphism $\varphi: B_{r,n}(c)\sra A$ defined by $e_i\mapsto e_i(p)$, 
maps the universal decomposition \eqref{eq1:undec} of $X^n(c)$ to $p(X)q(X)
$.
As a matter of 
$$
q(X)={X^n(c)\over p(X)}={\phi(X^n(c))\over \phi(\ttp_r(X))}=\phi\left(X^n(c)\over \ttp_r(X)\right)=\phi(X^{n-r}+h_1(c)X^{n-r-1}+\cdots+h_{n-r}(c)).
$$
\qed

A crucial result for our analysis is provided by the following:

\bclm{\bf Theorem.} \label{thm:lakth}{\em  There is a  $B_0(c)$-module isomorphism
$$
\bw^rV(c)\lra B_r(c),
$$
which turns $\bw^rV(c)$ into a free $B_r(c)$-module of rank $1$ generated by $\wX^r_0$.
}
\eclm
\proof 
It is a consequence of \cite[Theorem 0.1 (The Main Theorem)]{LakTh01}, once 
one rename $n$ by $r$ and identifies $B_0(c)$ with the algebra $S$ of 
symmetric functions in $r$-indeterminates over the ring $A=B_0(c)$. Notice 
that the above identification is possible due to the fundamental theorem of 
symmetric functions \cite[(2.4) and Remark following it]{MacDonald}\qed

\bclm{\bf Definition.}\label{def1:def17} {\em Let $g_0(X),g_1(X),
\ldots,g_{r-1}(X)$ elements of $V(c)[[X^{-1}]]$. The {\em residue} of the 
ordered $r$-tuple $(g_{r-1}(X),g_{r-2}(X),g_{0}(X))$ is defined as follows. 
If $r=1$,
	\begin{center}
		$\Res (g_0(X))=$ coefficient in $B_0(c)$ of the power $X^{-1}$ 
		occurring in $g_0(X)$.
\end{center}
For $r>1$:
\begin{eqnarray*}
\Res\left(g_{r-1}(X),\ldots, g_{0}(X)\right)=\begin{vmatrix}\Res(g_{r-1}
(X))&\ldots&\Res(g_0(X))\cr\cr
	\Res(Xg_{r-1}(X))&\ldots&\Res(Xg_0(X))\cr
	\vdots&\ddots&\vdots\cr
	\Res(X^{r-1}g_{r-1}(X))&\ldots&\Res(X^{r-1}g_0(X))\end{vmatrix}.
	\end{eqnarray*}
}
\eclm

Let now  $f_0(X),f_1(X),\ldots,f_{r-1}(X)\in V(c)$ and recall notation  \eqref{eq1:ttprX}.
\bclm{\bf Theorem.}\label{thm:resthm}
\be
f_{r-1}(X)\w\cdots\w f_0(X)=\Res\left({f_{r-1}(X)\over\ttp_r(X)},\ldots,{f_{0}(X)\over\ttp_r(X)}\right)\wX^r_0\label{eq:resfor}
\ee
\eclm
\proof \cite[Main Theorem]{LakTh01}.

\section{The $B_r(c)$-module structure of $\bw^rV(c)$}\label{sec2:sec2}

\claim{} \label{sec2:glnwr} The module structure alluded to in Theorem~\ref{thm:lakth} can be described as follows. 
Let us consider the standard representation of the Lie algebra 
$$
\delta:gl(V(c))\sra \End_{B_0(c)}(\bw V(c)),
$$ 
mapping $\phi\mapsto \delta(\phi)$, where $\delta(\phi)(u)=\phi(u)$ for all $u\in V(c)$ and
$$
\delta(\phi)(v\w w)=\delta(\phi)v\w w+ v\w \delta(\phi)w,\qquad \forall v,w\in\bw V(c). 
$$
In particular, we have:
\bclm{\bf Proposition.} {\em
The exterior power $\bw^rV(c)$ is an $gl(V(c))$-module.
}
\eclm
\proof Indeed, one gets a representation via the map $\phi\mapsto \delta(\phi)_{|\bw^rV(c))}$.\qed

\claim{} Relying on the above representation, we now describe more explicitly the module structure of $\bw^rV(c)$ over $B_r(c)$. Consider the endomorphism of $V(c)$ given by the multiplication by $X$. Define
$$
\sigma_+(z):=\sum_{i\geq 0}\sigma_iz^i=\exp\left(\sum_{i=1}{1\over i}\,\,\delta(X^i)z^i\right):\bw V(c)\lra \bw V(c)[[z]],
$$
which is the unique Hasse-Schmidt (HS) derivation on $\bw V(c)$ \cite{HSDGA} such that
$$
{\sigma_i}_{|V(c)}:V(c)\sra V(c),
$$
 is the multiplication by $X^i$. 
We denote by $\ovsig_+(z)$ the inverse HS derivation:
$$
\ovsig_-(z):=\sum_{i\geq 0}(-1)^i\ovsig_iz^i=\exp\left(-\sum_{i=1}{1\over i}\,\,\delta(X^i)z^i\right):\bw V(c)\lra \bw V(c)[[z]].
$$

Then, the $B_r$-module structure of $\bw^rV(c)$ is achieved by declaring that

\be
h_iu=\sigma_iu\qquad \mathrm{or, \,\, equivalently}\qquad e_iu=\ovsig_iu,\label{eq:IImosstr}
\ee
for all $u\in\bw^rV(c)$, where
 $h_i\in B_r$ is defined by \eqref{eq:Hrz}. This is precisely the same considered by Laksov and Thorup as one infer by comparing the references \cite[Theorem 2.4]{SCHSD}  and  \cite[Sections 2.1--2.2]{LakTh01}.
 
In the following result we prove the Giambelli's formula\cite[Formula 5.43]{HSDGA}
\bclm{\bf Theorem.}\label{thmgmblt} {\em Let 
\begin{eqnarray*}
\Delta_\blamb(H_r(c))=\begin{vmatrix}	
		h_{\lambda_1}(c) & h_{\lambda_2-1}(c) &\ldots & h_{\lambda_r -r+1}(c) \\
		h_{\lambda_1+1}(c) & h_{\lambda_2}(c) &\ldots & h_{\lambda_r-r+2}(c) \\
		\vdots  & \vdots & \ddots & \vdots \\
		h_{\lambda_1+r-1}(c) & h_{\lambda_2+r-2}(c) &\ldots &  h_{\lambda_r}(c) \\
	\end{vmatrix}.
\end{eqnarray*}
Then $\wXc^r_\blamb=\Delta_\blamb(H_r(c))\wXc^r_0$.}
\eclm
\proof  By Laksov--Thorup residue formula \cite[Theorem 0.1 (2)]{LakTh01}:
\begin{eqnarray*}
	\wXc^r_\blamb&=&X^{r-1+\lambda_1}(c)\w X^{r-2+\lambda_2}(c)\w\cdots\w X^{\lambda_r}(c)\cr\cr
	&=&\Res\left({X^{r-1+\lambda_1}(c)\over \ttp_r(X)},{X^{r-2+\lambda_2}(c)\over \ttp_r(X)}, \ldots,{X^{\lambda_r}(c)\over \ttp_r(X)}\right)\wX^r_0.
\end{eqnarray*}
Now:
\begin{eqnarray*}
	&&\Res \left({X^k \cdot X^{r-j+\lambda_j}(c)\over \ttp_r(X)}\right)\cr\cr&=&
	\Res \left({X^k \cdot X^{r-j+\lambda_j}(c)\over X^r}\left(1+{h_1 \over X}+{h_2 \over X^2} + \cdots \right)\right)\cr\cr
	&=&\Res  X^{k-r}\left( X^{r-j+\lambda_j}-c_1 X^{r-j+\lambda_j-1}+c_2 X^{r-j+\lambda_j-2}-\cdots 
	+(-1)^{r-j+\lambda_j}c_{r-j+\lambda_j}\right)\cr\cr
	&\cdot & 	\left(1+{h_1 \over X}+{h_2 \over X^2} + \cdots \right) \cr\cr
	&=&h_{\lambda_j+k-j+1}-c_1 h_{\lambda_j+k-j}+c_2 h_{\lambda_j+k-j-1}- \cdots + (-1)^{r-j+\lambda_j}c_{r-j+\lambda_j}\cr\cr
	&=&h_{\lambda_j+k-j+1}(c).
\end{eqnarray*}

The  definition \ref{def1:def17} of residue  and a simple computation gives that

\begin{eqnarray*}
	&&\Res\left({X^{r-1+\lambda_1}(c)\over \ttp_r(X)},{X^{r-2+\lambda_2}(c)\over \ttp_r(X)}, \ldots,{X^{\lambda_r}(c)\over \ttp_r(X)}\right)\hskip 20pt
	\cr\cr\cr &=&\begin{vmatrix}	
		h_{\lambda_1}(c) & h_{\lambda_2-1}(c) &\ldots & h_{\lambda_r -r+1}(c) \\
		h_{\lambda_1+1}(c) & h_{\lambda_2}(c) &\ldots & h_{\lambda_r-r+2}(c) \\
		\vdots  & \vdots & \ddots & \vdots \\
		h_{\lambda_1+r-1}(c) & h_{\lambda_2+r-2}(c) &\ldots &  h_{\lambda_r}(c) \\
	\end{vmatrix}=\Delta_\blamb(H_r(c)).
 \end{eqnarray*}
 \qed

\bclm{\bf Proposition.} {\em  The Schur determinants $\Delta_\blamb(H_r(c))$ are $B_0(c)$-bases of $B_r(c)$. More precisely \be
B_r(c)=\bigoplus_{\blamb\in\Pcal_r}B_0(c)\cdot\Delta_\blamb(H_r(c)).
\ee
}
\eclm
\proof Indeed, a $B_0(c)$-basis of $\bw^rV(c)$ is given either by $\wXc^r_\blamb$ or $\wXc^r_\blamb$. Thus if $P(e_1,\ldots,e_r)\in B_r(c)$, then 
$P(e_1,\ldots,e_r)\wXc^r_0\in \bw^rV(c)$ and can be then uniquely expressed as a  $B_0(c)$--linear combination of $\wXc^r_\blamb$ ($\wXc^r_\blamb$ respectively).
$$
P(e_1,\ldots,e_r)\wX^r_0=\sum_{\blamb}a_\blamb\wXc^r_\blamb=\sum_{\blamb}a_\blamb\Delta_\blamb(H_r)\wXc^r_0,
$$
from which 
$$
P(e_1,\ldots,e_r)=\sum_{\blamb}a_\blamb\Delta_\blamb(H_r),
$$
for some unique $a_\blamb\in B_0(c)$. \qed

\bclm{\bf Corollary.} {\em The ring $B_r(c)$ is a $gl(V(c))$-module}.
\eclm
\proof It follows from transporting the $gl(V(c))$-module structure via the isomorphism claimed in Theorem~\ref{thm:lakth}.\qed

\section{The $gl(V(c))$-structure of $B_r(c)$}\label{sec3:sec3}

The purpose of this section is to appropriately define the $gl(V(c))$-structure of $B_r(c)$, having in view the description of the structure of $B_{r,n}$ as a module over the Lie algebra $gl(V_n(c))$, where $V_n(c):=V(c)/(X^n(c))$. 
The isomorphism $B_r(c)\sra \bw^rV(c)$ implies the following:
\bclm{\bf Corollary.} {\em The ring $B_r(c)$ is a $gl(V(c))$-module.}
\eclm
\proof It follows from transporting the $gl(V(c))$-module structure of $\bw^rV(c)$ via the isomorphism claimed in Theorem~\ref{thm:lakth}.\qed

\bclm{\bf Definition.}\label{def:mndef}  {\em Let $f,g\in B_{0}(c)[t]$. The equality:
	\be
	\big[(f(X)\otimes g(\d))\star\Delta_\blamb(H_r(c))\big]\wXc^r_0=f(X)\w g(\delta)\lrcorner \wXc^r_\blamb,
	\ee
	defines $gl(V_n(c))$-structure of $B_r(c)$, i.e. the left hand side is defined through the right hand side.}

\claim{}\label{maintool} Our main tool to determine the explicit structure within a more general context than the one dealt with by Gatto and  Salehyan is a clever way to express the classical well known contraction of an exterior monomial against (a generating series of) linear forms. 
Let $\d\in V(c)^*$ arbitrarily chose. Then an easy application of the definition shows that
\be
\d\lrcorner ( P_1(X)\w \cdots\w P_r(X))=\begin{vmatrix}\d(P_1(X)&\cdots&\d(P_r(X))\cr\cr P_1(X)&\cdots&P_r(X)\end{vmatrix}\label{eq:conset}.
\ee

The right-hand side of \eqref{eq:conset} means  an alternating sum  where the coefficient of $\d(P_i(X))$ is the element of $\bw^{r-1}V(c)$ obtained by removing the $j$-th wedge factor from $P_1(X)\w\cdots\w P_r(X)$. For instance, the coefficient of $\d(P_1(X))$ is $P_2(X)\w\cdots\w P_r(x)$ while that of $\d(P_2(X)$ is $-P_1(X)\w P_3(X)\w\cdots\w P_r(X)$.
\bclm{\bf Example.}

\medskip
\begin{eqnarray*}
	\d \lrcorner (P_1(X)\w P_2(X)\w P_3(X))&=&\begin{vmatrix}\d(P_1(X))&\d(P_2(X))&\d(P_3(X))\cr\cr\cr
		P_1(X)&P_2(X)&P_3(X)\end{vmatrix}\cr\cr\cr&=&\d(P_1(X))\cdot P_2(X)\w P_3(X)-\d(P_2(X))\cdot P_1(X)\w P_3(X)\cr\cr &+&\d(P_3(X)\cdot P_1(X)\w P_2(X).
\end{eqnarray*}
Notice that $\d(P(X))\in B_0(c)$.
\eclm

\bclm{\bf Example.} Using this formalism, let us see how Definition~\ref{def:mndef} works by computing  $(X^3\otimes \d^2)\star\Delta_{(2,1)}(H_2(c))$.
By definition:
\begin{eqnarray*}
	(X^3\otimes \d^2)\star\Delta_{(2,1)}(H_2(c))\wXc^2_0&=&X^3\w (\d^2 \lrcorner (X^{3}(c)\w X^1(c))\cr\cr
	&=&X^3\w \begin{vmatrix}\d^2(X^3(c))&\d^2(X(c))\cr\cr X^3(c)&X(c)
	\end{vmatrix}
	=X^3\w \begin{vmatrix}-c_1&0\cr\cr X^3(c)&X(c)
	\end{vmatrix}
	\cr	\cr\cr
	&=&-c_1X^3\w X(c)
	=-c_1X^3\w X+c_1^2X^3\w 1\cr\cr
	&=&(-c_1\Delta_{(2,1)}(H_2)+c_1^2h_2)\wXc^2_0,
\end{eqnarray*}
from which, finally:
$$
(X^3\otimes \d^2)\star\Delta_{(2,1)}(H_2(c))=-c_1\Delta_{(2,1)}(H_2)+c_1^2h_2=-c_1(h_1h_2-h_3)+c_1^2h_2.
$$
\eclm
\claim{} Define $\Ecal(z,w^{-1}):B_r(c)\sra B_r(c)[[z,w^{-1}]$ through:
\be
(\Ecal(z,w)\star\Delta_\blamb(H_r(c)))\wX^r_0=\sum_{i,j\geq 0}X^iz^i\w\d^jw^{-j}\lrcorner \wXc^r_\blamb.\label{eq1:defac}
\ee
\eclm
Our purpose is to find an expression for $\Ecal(w,c)\Delta_\blamb(H_r(c))$ in such a way that $(X^i\otimes\d^j)\Delta_\blamb(H_r(c))$ is the coefficient of $z^iw^{-j}$ in the product $\Ecal(z,w)\Delta_\blamb(H_r(c))$.

\bclm{\bf Lemma.} {\em  Let $f\in V(c)$. Then
$$
{\bm\partial}(w^{-1})(f(X))=f(w^{-1}).
$$
}
\eclm
\proof
Assume $f$ of degree $d\geq 1$. Then $f(X)=\sum_{i=0}^da_iX^d$, where $a_i\in 
B_0(c)$. Thus
\begin{eqnarray*}
	{\bm\d}(w^{-1})(f(X))&=&\sum_{j\geq 0}w^{-j}\d^j(\sum_{i=0}^da_iX^i)\cr\cr
	&=&\sum_{i=0}^da_i\sum_{j\geq 0}w^{-j}\d^j(X^i)\cr\cr
	&=&\sum_{i=0}^da_i\sum_{j\geq 0}w^{-j}\delta^{ij}=\sum_{i=0}^da_iw^{-i}=f(w^{-1}).\hskip50pt \qed
\end{eqnarray*}
As a consequence, we get
\bclm{\bf Corollary.} ${\bm\partial}(w^{-1})(X^j(c))=w^{-j}(c)$.\qed
\eclm
\bclm{\bf Proposition.}\label{lem:lem26}{ \em
	\begin{eqnarray}
		&&\hskip-20pt {\bm \d}(w^{-1})\lrcorner (f_{r-1}(X)\w\cdots\w f_0(X))
		=\begin{vmatrix}f_{r-1}(w^{-1})& f_{r-2}(w^{-1})&\ldots&f_0(w^{-1})\cr\cr\cr
			f_{r-1}(X)& f_{r-2}(X)&\ldots&f_0(X)\cr\end{vmatrix}\cr\cr\cr
		&=&\begin{vmatrix}f_{r-1}(w^{-1})& f_{r-2}(w^{-1})&\ldots&f_0(w^{-1})\cr\cr
			\Res\,\displaystyle{f_{r-1}(X)\over \ttp_{r-1}(X)}&\Res\displaystyle{f_{r-2}(X)\over  \ttp_{r-1}(X)}&\ldots&\Res\displaystyle{f_{0}(X)\over \ttp_{r-1}(X)}\cr\cr
			\vdots&\vdots&\ddots&\vdots\cr\\
			\Res\displaystyle{X^{r-2}\cdot f_{r-1}(X)\over \ttp_{r-1}(X)}&\Res\displaystyle{X^{r-2}\cdot f_{r-2}(X)\over \ttp_{r-1}(X)}&\ldots&\Res\displaystyle{X^{r-2}f_{0}(X)\over \ttp_{r-1}(X)}\end{vmatrix}\wXc^{r-1}_0.\hskip40pt \label{eq:lakgmb}
	\end{eqnarray}
}\eclm

\proof The first equality of \eqref{eq:lakgmb} is obvious by virtue of \eqref{eq:conset}. To prove the second, first observe that $f_{r-1}(w^{{-1}})$ is the coefficient of the determinant
\begin{eqnarray*}
	&&\begin{vmatrix}\Res\,\displaystyle{f_{r-2}(X)\over  \ttp_{r-1}(X)}&\ldots& \Res\,\displaystyle{f_{0}(X)\over  \ttp_{r-1}(X)}\cr\vdots&\ddots&\vdots\cr\cr \Res\,\displaystyle{X^{r-2}f_{r-2}(X)\over  \ttp_{r-1}(X)})&\ldots& \Res\,\displaystyle{X^{r-2}f_{0}(X)\over  X^{r-2}\ttp_{r-1}(X)}
	\end{vmatrix}\wXc^{r-1}_0=\cr\cr\cr
	&=&\Res\left(f_{r-2}(X),\ldots, f_{0}(X)\right)\wXc^{r-1}_0=f_{r-2}(X)\w\cdots\w f_{0}(X),
\end{eqnarray*}
which is precisely the coefficients of $f_{r-1}(w^{-1})$ in the first member of \eqref{eq:lakgmb}. Then the property is proved because the determinant is alternating with respect to the columns.\qed


\bclm{\bf Definition.} {\em
	For all $\blamb\in \Pcal_r$, define 
	\be
	\Delta_\blamb(\bfw^{-1}(c),H_{r-1}(c)):={w^{r-1}}\,{\bm\d}(w^{-1})\star \Delta_\blamb(H_r).\label{eq:3411}
	\ee
}
\eclm

\bclm{\bf Corollary.} {\em The explicit expression of the right hand side of \eqref{eq:3411} is:
	\begin{eqnarray*}
		\Delta_\blamb(\bfw^{-1}(c),H_{r-1}(c))
		&=&{w^{r-1}}\begin{vmatrix}w^{-r+1-\lambda_1}(c)&w^{-r+2-\lambda_2}(c)&\cdots& w^{-\lambda_r}(c)\cr\cr
			h_{\lambda_1+1}(c)&h_{\lambda_2}(c)&\cdots&h_{\lambda_r+r-2}(c)\cr\cr\vdots&\vdots&\ddots&\vdots
			\cr\cr
			h_{\lambda_1+r-1}(c)&h_{\lambda_2+r-2}(c)&\cdots&h_{\lambda_r}(c)
		\end{vmatrix}.
	\end{eqnarray*}
}
\eclm

\proof It is a straightforward consequence of Proposition~\ref{lem:lem26}, by choosing $f_{r-j}(X)=X^{r-j+\lambda_j}(c)$, for all $1\leq j\leq r$. \qed

For the next proposition it is useful to remark $\wX^k_0=\wXc^k_0$.

\bclm{\bf Lemma.} \label{lemw} {\em The following equality holds for all $k\geq 1$ and all $\blamb\in\Pcal_r$:
	\begin{eqnarray*}
		&&\wX^{k}_0\w \begin{vmatrix}f_{r-1}(w^{-1})&f_{r-2}(w^{-1})&\cdots& f_0(w^{-1})\cr\cr
			f_{r-1}(X)&f_{r-2}(X)&\ldots&f_0(X)\end{vmatrix}\cr\cr\cr
		&=&(-1)^{r+k-1}\begin{vmatrix}f_{r-1}(w^{-1})&f_{r-2}(w^{-1})&\cdots& f_0(w^{-1})&0&\ldots&0\cr\cr
			f_{r-1}(X)&f_{r-2}(X)&\ldots&f_0(X)&X^{k-1}&\ldots&X^0
		\end{vmatrix}.\\ 
	\end{eqnarray*}
	
}
\eclm
\proof By just comparing the coefficients of the expression $f_j(w^{-1})$ ($0\leq j\leq r-1$) on eithe side of the equality. \qed

\bclm{\bf Corollary.}
\begin{eqnarray*}
	&&\wX^{k}_0\w \begin{vmatrix}w^{-r+1-\lambda_1}(c)&w^{-r+2-\lambda_2}(c)&\cdots& w^{-\lambda_r}(c)\cr\cr
		X^{r-1+\lambda_1}(c)&X^{r-2+\lambda_2}(c)&\cdots&X^{\lambda_r}(c)
	\end{vmatrix}\cr\cr\cr
	&=&(-1)^{r+k-1}\begin{vmatrix}w^{-r+1-\lambda_1}(c)&w^{-r+2-\lambda_2}(c)&\cdots& w^{-\lambda_r}(c)&0&\ldots&0\cr\cr
		X^{r-1+\lambda_1}(c)&X^{r-2+\lambda_2}(c)&\cdots&X^{\lambda_r}(c)&X^{k-1}&\ldots&X^0
	\end{vmatrix}\\ &=&\wXc^{r+k}_{\blamb-(k^r)}.
\end{eqnarray*}

\eclm

\proof The formula is a straightforward consequence of Lemma~\ref{lemw}. Once one checks the formula for $k=1$, the other cases follow by induction. Clearly the result is zero if $\blamb-(k^r)$ were not a partition, i.e. if $\lambda_r-k<0$. \qed

\bclm{\bf Proposition.}\label{prop:lem} {
	\begin{eqnarray*}
		&&\sum_{i\geq 0} X^iz^i\w {\bm \d}(w^{-1}(c))\lrcorner \wXc^r_\blamb\cr&=&{z^{r-1}\over w^{r-1}}{1\over E_r(z)}\begin{vmatrix}w^{-\lambda_1}(c)&w^{1-\lambda_2(c)}&\ldots&w^{r-1-\lambda_r}(c)&0\cr\cr
			\ovsig_-(z)\cdot X^{r+\lambda_1}(c)&\ovsig_-(z) X^{r-1+\lambda_2}(c)&\ldots&\ovsig_-(z) X^{\lambda_r+1}(c)&1\end{vmatrix}\cr\cr\cr\cr
		&=&{z^{r-1}\over w^{r-1}}{1\over E_r(z)}\begin{vmatrix}w^{-\lambda_1}&\ldots&w^{r-1-\lambda_r}&0\cr\cr X^{r+\lambda_1}(c)-\displaystyle{X^{r-1+\lambda_1}(c)\over z}&\ldots&  X^{1+\lambda_r}(c)
			-\displaystyle{X^{\lambda_r}(c)\over z}&1\end{vmatrix}.
	\end{eqnarray*}
}
\eclm
\proof
It follows by noticing that

\begin{eqnarray*}
	&&\sum_{i\geq 0} X^iz^i\w {\bm \d}(w^{-1})\lrcorner \wXc^r_\blamb
	\cr\cr
	&=&\sigma_+(z)X^0\w w^{r-1}\begin{vmatrix}w^{-\lambda_1}&w^{1-\lambda_2}&\ldots&w^{r-1-\lambda_r}&0\cr\cr  X^{r-1+\lambda_1}(c)&X^{r-2+\lambda_2}(c)&\ldots&X^{\lambda_r}(c)&1\end{vmatrix}\cr\cr\cr\cr
	&=&w^{r-1}\sigma_+(z)\left(X^0\w \ovsig_+(z)\begin{vmatrix}w^{-\lambda_1}&w^{1-\lambda_2}&\ldots&w^{r-1-\lambda_r}\cr\cr  X^{r-1+\lambda_1}(c)&X^{r-2+\lambda_2}(c)&\ldots&X^{\lambda_r}(c)\end{vmatrix}\right),
	%
	%
\end{eqnarray*}
where in the last equality we have used integration by parts. The fact that $\ovsig_+(z)$ is a HS derivation translate in filtering the $\ovsig_+(z)$ inside the diagram, by letting it acting singularly on each element of the second row:

\begin{eqnarray*}
	&=&w^{r-1}\sigma_+(z)\left(X^0\w\begin{vmatrix}w^{-\lambda_1}&w^{1-\lambda_2}&\ldots&w^{r-1-\lambda_r}\cr\cr   \ovsig_+(z)X^{r-1+\lambda_1}(c)& \ovsig_+(z)X^{r-2+\lambda_2}(c)&\ldots& \ovsig_+(z)X^{\lambda_r}(c)\end{vmatrix}\right).
\end{eqnarray*}
The last equality is equivalent to
\begin{eqnarray}
	&=&w^{-r+1}z^{r-1}\sigma_+(z)\begin{vmatrix}w^{-\lambda_1}&w^{1-
	\lambda_2}&\ldots&w^{r-1-\lambda_r}&0\cr\cr   \ovsig_-(z)X^{r+\lambda_1}
	(c)& \ovsig_-(z)X^{r-1+\lambda_2}(c)&\ldots& \ovsig_-(z)X^{1+\lambda_r}
	(c)&1\end{vmatrix}\cr\cr\cr
	&=&{z^{r-1}\over w^{r-1}}\sigma_+(z)\ovsig_-(z)\begin{vmatrix}w^{-
	\lambda_1}&w^{1-\lambda_2}&\ldots&w^{r-1-\lambda_r}&0\cr\cr   X^{r+
	\lambda_1}(c)& X^{r-1+\lambda_2}(c)&\ldots&X^{1+\lambda_r}
	(c)&1\end{vmatrix}.
	\label{eq:psdet}
\end{eqnarray}
Since the ``determinant'' occurring in \eqref{eq:psdet} is a linear 
combination of vectors of $\bw^{r}V$, which is an eigenspace of $\sigma_+(z)
$ having $\displaystyle{1\over E_r(z)}$ as eigenvalue, we eventually have:
$$
{z^{r-1}\over w^{r-1}}{1\over E_r(z)}\ovsig_-(z)\begin{vmatrix}w^{-
\lambda_1}&w^{1-\lambda_2}&\ldots&w^{r-1-\lambda_r}&0\cr\cr   X^{r+
\lambda_1}(c)& X^{r-1+\lambda_2}(c)&\ldots&X^{1+\lambda_r}
(c)&1\end{vmatrix}.\hskip100pt \qed
$$

\bclm{\bf Theorem.}\label{thm.thm210} {\em One has:
	\be
	\Ecal(z,w^{-1})\Delta_\blamb(H_r(c))= {z^{r-1}\over w^{r-1}}{1\over 
	E_r(z)}\Delta_\blamb(\bfw^{-1}(c),\ovsig_-(z)H_r(c)).
	\label{eq:finalformula}
	\ee
}\eclm
\proof
By Definition~\ref{def:mndef}, the formal power series $\Ecal(z,w^{-1})
\Delta_\blamb(H_r(c))$ is the unique element of $B_r[[ z,w^{-1}]]$ such that
$$
\Ecal(z,w^{-1})\Delta_\blamb(H_r(c))\wXc^r_0=\sum_{i\geq 0} X^iz^i\w {\bm 
\d}(w^{-1})\lrcorner \wXc^r_\blamb.
$$
Now, by Propostion~\ref{prop:lem}, one has
\begin{eqnarray*}
	&&\sum_{i\geq 0}X^iz^i\otimes {\bm\d}(w^{-1})\lrcorner \wXc^r_\blamb\cr\cr\cr
	&=&{z^{r-1}\over w^{r-1}}{1\over E_r(z)}\ovsig_-(z)\begin{vmatrix}w^{-\lambda_1}&w^{1-\lambda_2}&\ldots&w^{r-1-\lambda_r}&0\cr\cr   X^{r+\lambda_1}(c)& X^{r-1+\lambda_2}(c)&\ldots&X^{1+\lambda_r}(c)&1\end{vmatrix},\cr\cr
\end{eqnarray*}
and by invoking Lemma~\ref{lem:lem26}, eventually:
\be
{z^{r-1}\over w^{r-1}}{1\over E_r(z)}\Delta_\blamb(\bfw^{-1}(c),\ovsig_-(z)H_r(c)),
\ee
(Notice that we are now with coefficients in $B_r(c)$). Now by ~\cite[Theorem 5.7]{pluckercone}, $\ovsig_-(z)$ can filter inside the determinant, to give finally the required expression~(\ref{eq:finalformula}).\qed
\claim{\bf Remark.} Notice that for $c_1=\ldots=c_k=0$, we obtain the main formula of \cite{gln}. However our proof is more general and neater and we will see how in a while will allow us to describe the $gl(V_n(c))$ structure of the universal decomposition algebra $B_{k,n}(c)$ of the polynomial $X_n(c)\in B_0(c)[X]$.

\bclm{\bf Example.} Let us check the formula for $r=3$, applied to $1\in B_r$. It gives:
\begin{eqnarray}
	&&\Ecal(z,w^{-1})1={z^{2}\over w^{2}}{1\over E_3(z)}\Delta_{(0)}(\bfw^{-1}(c),H_3(c))\cr\cr\cr
	&=&{z^2\over E_3(z)}\begin{vmatrix}w^{-2}(c)&w^{-1}(c)&1\cr\cr h_1(c)-\displaystyle{1\over z}&1&0\cr\cr h_2(c)-\displaystyle{h_1(c)\over z}&h_1(c)-\displaystyle{1\over z}&1\end{vmatrix}\cr\cr\cr
	&=&{z^2\over E_3(z)}\left\{w^{-2}(c)-w^{-1}(c)\left(h_1(c)-{1\over z}\right)+\left[\left(h_1(c)-{1\over z}\right)^2-h_2(c)+{h_1(c)\over z}\right]\right\}\cr\cr\cr
	&=&{1\over E_3(z)}\left(1-h_1z+(h_1^2-h_2)z^2+{z\over w}(1-h_1z)+{z^2\over w^2}\right)\cr\cr\cr
	&=&{1\over E_3(z)}\left(E_2(z)+{z\over w}E_1(z)+{z^2\over w^2}\right).\label{eq:allbm}
	\end{eqnarray}
\eclm
For example 
\begin{eqnarray*}
((X^5\otimes \d^1)1)\wX^3_0&=&X^5\w (\d^1\lrcorner X^2\w X^1\w X^0)\cr\cr
&=&-X^5\w X^2\w X^0=-\Delta_{(3,1)}(H_3)\wX^3_0,
\end{eqnarray*}
 from which
$(X^5\otimes \d^1)1=-\Delta_{(3,1)}(H_3)=h_4-h_1h_3$, which in turns is  precisely  the coefficient of $z^5w^{-1}$ in formula \eqref{eq:allbm}.

\claim{} Let $\bfs(w)=\displaystyle{\sum_{i\geq 0}}s_iw^i:=\displaystyle{1\over \bfc(w)}$
and define  $(\d^j(s))_{j\geq 0}$  through the equality:
\be
\sum_{i\geq 0}\d^j(s)w^{-j}=\bfs(w)\d(w^{-1})={\d(w^{-1})\over \bfc(w)}.\label{eq:eq16}
\ee
Then $(\d^j(s))_{j\geq 0}$ is the dual basis of $\bfX(c)=(X^i(c))_{i\geq 0}$ as it is easy to check.
\bclm{\bf Corollary.} \label{cor3:cor22} {\em  In terms of the adapted bases $\wXc$ and $(\d^j(s))_{j\geq 0}$ one has:
\be
(\sum_{j\geq 0}X^j(c)z^j\otimes \d^j(s)w^{-j})\star\Delta_\blamb(H_r(c))={z^{r-1}\over w^{r-1}}{c(z)\over c(w)}{1\over E_r(z)} \Delta_\blamb(\bfw^{-1},\ovsig_-(z)H_r(c)).\label{eq:eqcrl}
\ee
}
\eclm
\proof
Indeed, $\displaystyle{\sum_{j\geq 0}}X^j(c)z^j=\bfc(z)\sigma_+(z)1$ and by \eqref{eq:eq16}:
\begin{eqnarray*}
	(\sum_{j\geq 0}X^j(c)z^j\otimes \d^j(s)w^{-j})\star\Delta_\blamb(H_r(c))&=&\left(\bfc(z)\sigma_+(z)1\otimes {\d^j(w^{-j})\over \bfc(w)}\right)\Delta_\blamb(H_r(c))\cr\cr
	&=&{\bfc(z)\over \bfc(w)}(\sigma_+(z)\otimes \d^j(w^{-1})\Delta_\blamb(H_r(c)),
\end{eqnarray*}
whence \eqref{eq:eqcrl}, after applying Theorem~\ref{thm.thm210}.\qed

\section{The $gl(V_n(c))$-structure of $B_{r,n}(c)$.}\label{sec4:sec4}
\claim{} For all $k\geq 0$, denote
$$
V_k(c)=\bigoplus_{j=0}^{k-1}B_0(c)\cdot X^j(c)\qquad \mathrm{and}\qquad V^k(c)=\bigoplus_{j\geq k}B_0(c)\cdot X^{j}(c).
$$
Therefore:
$$
V(c)=V_k(c)\oplus V^k(c).
$$

By  \cite[Main Theorem 0.1]{LakTh01}, or our previous remarks,  there is a $B_0(c)$-module isomorphism between the universal decomposition algebra 
$$
B_{r,n}(c):={B_r(c)\over I_{r,n}(c)},
$$ 
of the monic polynomial $X^n(c)$, where
$I_{r,n}(c)$ is the ideal $(h_{n-r+1}(c),\ldots, h_n(c))$, and  $\bw^rV_n(c)$, where
$$
V_n(c):={V(c)\over (X^n(c))}={B_0(c)[X]\over (X^n(c))}=\bigoplus_{i=0}^{n-1}B_0(c)\cdot X^i(c).
$$
Since $\bw^rV_n(c)$ represents $gl(V_n(c))$ via the map $\delta:gl(V_n(c))\sra End_{B_0(c)}\bw V_n(c))$, it turns out that $B_{r,n}(c)$ is itself a $gl(V_n(c))$-module.
Let $\pi_{r,n}:B_{r}(c)\sra B_{r,n}(c)$ the canonical epimorphism.
\bclm{\bf Lemma.} \label{lem:lem352} {\em Suppose $\blamb\in\Pcal_r$ such that $\lambda_1\geq n-r+1$. Then $\Delta_\blamb(H_r(c))\in I_{r,n}(c)$.
}
\eclm

\proof It is  an adaptation of similar proof in the book \cite[p.~113--114]{HSDGA}.
Recall that by its very definition,  $\Delta_\blamb(H_r(c))$ belongs to the ideal $(h_{\lambda_1}(c), h_{\lambda_1+1}(c),\ldots,h_{\lambda_1+r-1}(c))$. Thus, if $\lambda_1\geq n-r+1$, $\Delta_\blamb(H_r(c))\in I_{r,n}(c)$, as claimed.\qed

Let
\be
\Delta_\blamb(H_{r,n}(c))=\Delta_\blamb(H_{r}(c))+I_{r,n}(c)\label{eq:defhrn}.
\ee
Clearly  $\Delta_\blamb(H_{r,n}(c))=0$ in $B_{r,n}(c)$ if $\lambda_1\geq n-r+1$. But if $\blamb\in\Pcal_{r,n}$, i.e. the Young diagram of $\blamb$ is contained in a $r\times (n-r)$ rectangle, then $\Delta_\blamb(H_{r,n}(c)$ are linearly independent, because  they are in $B_r$. So $(\Delta_\blamb(H_{r,n}(c))$ is a $B_0(c)$-basis of $B_{r,n}(c)$.

We consider
the generating function 
\be
\Ecal_n(z,w^{-1})=\sum_{i,j= 0}^{n-1}X^i(c)z^i\otimes \partial^j(s)w^{-j}:B_{r,n}\sra B_{r,n}[z,w^{-1}].
\ee Clearly this is a polynomial in $z$ and $w^{-1}$, because $(X^0,X(c),\ldots,X^{n-1}(c))$ is a basis of $V_n(c)$ and $X^{n+i}(c)=0$ for all $i\geq 0$. 

\bclm{\bf Proposition.}\label{prop:prop43} {\em For all $0\leq j\leq n-1$
	$$
	X^i(c)\otimes \d^j(s)(I_{r,n}(c))\subseteq I_{r,n}(c).
	$$
}
\eclm
\proof
First recall that $I_{r,n}(c)\bw^rV(c)= V^n(c)\w \bw^{r-1}V(c)$. In fact
for all $\blamb\in\Pcal_r$ one has
$$
I_{r,n}(c)\wXc^r_\blamb=I_{r,n}(c)\Delta_\blamb(H_r(c)\wXc^r_0\subseteq I_{r,n}(c)\wXc^r_0. 
$$
So, it suffices to show that $I_{r,n}(c)\wXc^r_0\in V^n(c)\w \bw^{r-1}V(c)$, i.e. that $h_{n-r+j}(c)\wXc^r_0\in V^n(c)\w \bw^{r-1}V(c)$ for all $1\leq j\leq r$.
But
$$
h_{n-r+j}(c)\wXc^r_0=X^{n+j-1}(c)\w X^{r-2}(c)\w\cdots\w X^0(c)\in V^n(c)\w\bw^{r-1}V(c),
$$
 which proves the inclusion $I_{r,n}(c)\bw^rV(c)\subseteq V^n(c)\w \bw^{r-1}V(c)$.  Conversely if $\wXc^r_\blamb\in V^n(c)\w\bw^{r-1}V(c)$, then $\lambda_1\geq n-r+1$. Thus $\wXc^r_\blamb=\Delta_\blamb(H_r(c))\wXc^r_0$, with $\Delta_\blamb(H_r(c))\in I_{r,n}(c)$, because $\lambda_1\geq n-r+1$ and Lemma~\ref{lem:lem352}.
Now if $j\leq n-1$
\begin{eqnarray*}
	(\d^j(s)\star I_{r,n}(c))\wXc^r_0&=&\d^j(s)(V^n(c)\w\bw^{r-1}V(c)\subseteq V^n(c)\w \d^j(s)\lrcorner \bw^{r-1}V(c)\cr\cr
	&\subseteq& V^n(c)\w\bw^{r-1}V(c)=I_{r,n}(c)\wXc^r,
\end{eqnarray*}
whence $\d^j(s)\star I_{r,n}(c)$. It follows that
\begin{eqnarray*}
	X^i(c)\otimes \d^j(s)\star(\Delta_\blamb(H_r(c))+I_{r,n})&=& X^i(c)\otimes \d^j(s)(\Delta_\blamb(H_r(c))+ I_{r,n}\cr\cr
	&=&V^n(c)\w \bw^{r-1}V(c). \hskip 60pt \qed
\end{eqnarray*}

We have then the following:
\bclm{\bf Theorem.}\label{thm44} {\em 
	Let $\Ecal_n(z,w^{-1})=\displaystyle{\sum_{i,j=0}^{n-1}}X^i(c)\otimes \partial_j(s)z^iw^{-j}$.
	Then
\be
	\Ecal_n(z,w^{-1})(\Delta_\blamb(H_{r,n}(c)))={z^{r-1}\over w^{r-1}}\left(1+\sum_{i=1}^{n-1}s_iw^i\right){c(z)\over E_r(z)}\pi_{r,n}\Delta_\blamb(\bfw^{-1}(c),\ovsig_-(z)H_r(c)).\label{eq:brn}
\ee
}
\eclm
\proof By definition of $\Delta_\blamb(H_{r,n}(c))$, one has:

\begin{center}
\begin{tabular}{lrll}
&&$\Ecal_n(z,w^{-1})(\Delta_\blamb(H_{r,n}(c)))$&\cr\cr
&$=$&$\Ecal_n(z,w^{-1})(\Delta_\blamb(H_r(c))+I_{r,n}(c))$&(Definition \eqref{eq:defhrn}\cr
&&&of $\Delta_\blamb(H_{r,n})$\cr\cr
&$=$&$\Ecal_n(z,w^{-1})(\Delta_\blamb(H_r(c)))+I_{r,n}(c)$&(Proposition~\ref{prop:prop43})\cr\cr
&$\hskip-5pt =$&$\hskip-10pt \displaystyle{z^{r-1}\over w^{r-1}}\left(\displaystyle{1\over c(w)}\,\mod\, w^{n}\right)\pi_{r,n}\displaystyle{c(z)\over E_r(z)}\Delta_\blamb(\bfw^{-1}(c),\ovsig_-(z)H_r(c)),$&(Corollary 3.12).\cr\cr
&$=$&$\hskip-10pt \displaystyle{z^{r-1}\over w^{r-1}}\left(\displaystyle{1\over c(w)}\,\mod\, w^{n}\right)\pi_{r,n}\displaystyle{c(z)\over E_r(z)}\pi_{r,n}\Delta_\blamb(\bfw^{-1}(c),\ovsig_-(z)H_r(c)),$&
\end{tabular}
\end{center}
where in the last equality we used the fact that $\pi_{r,n}$, the canonical 
epimorphism, is a ring homomorphism.\qed
\bclm{\bf Example.} Let us check the formula for $n=4$ and $r=2$ by 
computing the generating function of the images of $\Delta_{(2,1)}(H_{2,4}
(c))$.
 Applying Theorem \ref{thm44}, we get
\begin{eqnarray}
	&&\Ecal_4(z,w^{-1})(\Delta_{(2,1)}(H_{2,4}(c)))=
	\Ecal_4(z,w^{-1})\Delta_{(2,1)}(H_2(c))+I_{2,4}(c)\cr\cr\cr
	&=&{z\over w}{c(z)\over c(w)}{1\over E_2(z)}\Delta_{(2,1)}(\bfw^{-1}(c),
	\ovsig_-(z)H_{2,4}(c))\cr\cr\cr
	&=&{z}{1\over c(w)}{c(z)\over E_2(z)}
	\begin{vmatrix}w^{-3}(c)&w^{-1}(c)\cr\cr h_3(c)-\displaystyle{h_2(c)
	\over z}&h_1(c)-\displaystyle{1\over z}\end{vmatrix}+I_{2,4}(c)\cr\cr\cr
	&=&{z}(1+s_1w+s_2w^2+s_3w^3)\left(1+h_1(c)z+ h_2(c)z^2 \right)\cdot\cr
	\cr
	&\cdot& \left(-c_3h_1(c) +{1\over z}(c_3-c_1h_2(c))+{1\over w}c_2h_1(c)\right.\cr\cr&&\left. -{1\over w^2}\,c_1h_1(c)+{1\over zw}(h_2(c)-c_2)+{h_1(c)\over w^3}+{c_1\over w^2 z}-{1\over w^3 z}\right)\cr\cr 
	&=&
(1+h_1(c)z+h_2(c)z^2)\left({1\over w^3}(h_1(c)z-1)+{h_2(c)\over w}\right)\cr\cr\cr
&=& {h_2(c)\over w} + h_1(c)h_2(c){z\over w}+ h_2(c)^2 {z^2\over w} -{1\over w^3} + (h_1(c)^2-h_2(c)){z^2\over w^3}\cr\cr
&&
   + h_1(c) h_2(c){z^3\over w^3}\label{eq4:eq25}\\ \cr\cr
   &=&h_2(c)\,{1\over w}+\Delta_{(2,1)}(H_{2,4}(c)){z\over w}+\Delta_{(2,2)}(H_{2,4}(c)){z^2\over w}\cr\cr
   &-&{1\over w^3}+\Delta_{(1,1)}(H_{2,4}(c)){z^2\over w^3}+ \Delta_{(2,1)}(H_{2,4}(c)){z^3\over w^3}.
\label{eq4:eq26}
\end{eqnarray}

\noindent
For instance a direct computation gives:
\begin{eqnarray*}
((X^2(c)\otimes \d^{-1}(s))\Delta_{(2,1)}(H_{2,4}(c))\wXc^2_0&=&X^2(c)\w (\d^{-1}(s)\lrcorner X^{3}(c)\w X^1(c))\cr\cr
&=&X^3(c)\w X^2(c)=X^3\w X^2(c)=h_2(c)^2,
\end{eqnarray*}
which is precisely the coefficient of $z^2w^{-1}$ in either \eqref{eq4:eq26} or \eqref{eq4:eq25}.
\eclm

\bibliographystyle{amsplain}

\medskip
\medskip

\parbox[t]{3in}{{\rm Ommolbanin~Behzad}\\ 
\smallskip \vspace{-15pt}

	{\tt \href{mailto:behzad@iasbs.ac.ir}{behzad@iasbs.ac.ir}}\\
	{\it Department of Mathematics\\
	Institute for Advanced Studies\\ in Basic Sciences (IASBS)\\ P.O.Box 45195-1159  Zanjan\\ IRAN}}
\parbox[t]{3in}{{\rm Abbas Nasrollah Nejad}\\ 

\smallskip \vspace{-15pt}
	{\tt \href{mailto:	abbasnn@iasbs.ac.ir
		}{	abbasnn@iasbs.ac.ir
	}}\\
	{\it Department of Mathematics\\
			Institute for Advanced Studies\\ in Basic Sciences (IASBS)\\ P.O.Box 45195-1159  Zanjan\\ IRAN}}

\end{document}